\documentclass{lms}

\usepackage{amsmath}
\usepackage{amsfonts}
\usepackage{amssymb}
\usepackage{longtable}

\newtheorem{theorem}{Theorem}
\newtheorem{prop}{Proposition}

\newtheorem{lem}{Lemma}
\newtheorem{tab}{Table}

\def\zr{\ltimes}

\def\Real{\mathbb{R}}
\def\Co{\mathbb{C}}

\def\g{\mathfrak{g}}
\def\h{\mathfrak{h}}

\def\so{\mathfrak{so}}
\def\spin{\mathfrak{spin}}

\def\sl{\mathfrak{sl}}
\def\gl{\mathfrak{gl}}
\def\su{\mathfrak{su}}
\def\u{\mathfrak{u}}
\def\sp{\mathfrak{sp}}
\def\f{\mathfrak{f}}
\def\t{\mathfrak{t}}
\def\z{\mathfrak{z}}
\def\k{\mathfrak{k}}
\def\R{\mathcal{R}}

\def\N{\mathcal{N}}
\def\C{\mathcal{C}}

\def\id{\mathop\text{\rm id}\nolimits}

\def\tr{\mathop\text{\rm tr}\nolimits}

\def\pr{\mathop\text{\rm pr}\nolimits}

\def\Ric{\mathop{{\rm Ric}}\nolimits}

\def\be{\begin{equation}}
\def\ee{\end{equation}}

\begin{document}

\title[Holonomy of Einstein pseudo-Riemannian manifolds]{Holonomy algebras of Einstein pseudo-Riemannian manifolds}
\author{Anton S. Galaev}

\classno{53C29; 53C50}

\maketitle

\begin{abstract}
The holonomy algebras of  Einstein not Ricci-flat pseudo-Riemannian
manifolds of arbitrary signature are classified. As  illustrating
examples, the cases of Lorentzian manifolds, pseudo-Riemannian
manifolds of signature $(2,n)$ and the para-quaternionic-K\"ahlerian
manifolds with non-zero scalar curvature are considered. Einstein
not Ricci-flat metrics of signature $(2,n)$  with all possible
holonomy algebras are given.

\vskip0.5cm

{\bf Keywords:} holonomy algebra; Einstein pseudo-Riemannian
manifold; para-quaternionic-K\"ahlerian manifold.
\end{abstract}

\section{Introduction}

The holonomy group of a pseudo-Riemannian manifold gives rich
information about the geometry of the manifold. This motivates the
classification problem for holonomy groups. For simplicity one
usually restricts the attention to the  connected component of the
holonomy group, then it is enough to consider the corresponding Lie
algebra, called the holonomy algebra. Even in this case the
classification problem in arbitrary signature seems to be
unsolvable. The complete solution exists only in the Riemannian case
\cite{Ber,Al1,Besse,Joyce07,Bryant2} and in the Lorentzian case
\cite{BB-I,Leistner,Galmetr}. For an arbitrary signature, only the
classification of irreducible holonomy algebras is known
\cite{Ber,Bryant2}. The general case cannot be reduced to the
irreducible one unless the metric is  of the Riemannian signature.
For pseudo-Riemannian manifolds of signature different from
Riemannian and Lorentzian ones, only some partial case are
considered \cite{BI97,BBnn,Quat,ind4,F-K,IRMA,GalDB,I99}. There is
also a classification of connected irreducible holonomy groups of
torsion-free affine connections \cite{M-Sch99}; the groups
corresponding to the Ricci-flat case are found in \cite{Arms}.

In this paper we restrict our attention to the holonomy algebras of
Einstein not Ricci-flat pseudo-Riemannian manifolds. By now the
results have been known only in the case of irreducible holonomy
algebras \cite{Bryant2} and in the Lorentzian signature
\cite{GalEin}.

We show that the holonomy algebra of an Einstein not Ricci-flat
pseudo-Riemannian manifold contains a significant reductive
subalgebra, which is again the holonomy algebra of an Einstein
pseudo-Riemannian manifold and which is known. This allows us to get
a complete description of the holonomy algebras of such manifolds in
arbitrary signature. For that it is enough to combine some modules
of the reductive part of the holonomy algebra. As the illustration,
we classify holonomy algebras of Einstein not Ricci-flat
pseudo-Riemannian manifolds of signature $(2,n)$ as well as of
para-quaternionic-K\"ahlerian manifolds with non-zero scalar
curvature. In signature $(2,n)$, we construct  examples of Einstein
not Ricci-flat metrics with all possible holonomy algebras. We also
discuss the construction of  Einstein metrics in other signatures;
this construction is technical and it should be given elsewhere.

\section{Background}

The theory of holonomy algebras of pseudo-Riemannian manifolds can
be found e.g. in \cite{Besse,Bryant2,Joyce07}. Here we collect
some facts that  will be used below.

Let $(M,g)$ be a connected pseudo-Riemannian manifold of signature
$(p,q)$. The holonomy group of $(M,g)$ at a point $x\in M$ is a
Lie group that consists of the pseudo-orthogonal transformations
given by the parallel displacements along piece-wise smooth loops
at the point $x$, and it can be identified with a Lie subgroup of
the pseudo-orthogonal Lie group ${\rm O}(p,q)$. The corresponding
Lie subalgebra of $\so(p,q)$ is called the holonomy algebras. If
the manifold $M$ is simply connected, then the holonomy group is
connected and it is uniquely defined by the holonomy algebra.

Let $\g\subset\so(p,q)$ be a subalgebra. The space of curvature
tensors of type $\g$ is defined as follows
$$\R(\g)=\left\{R\in\wedge^2 (\Real^{p,q})^*\otimes
\g\,\left|\begin{array}{c}R(X, Y)Z+R(Y, Z)X+R(Z, X)Y=0\\ \text{
}\text{ for all  } X,Y,Z\in \Real^{p,q}\end{array}\right\}\right..$$
The above identity is called {\emph the first Bianchi identity}. Let
$L(\R(\g))\subset\g$ be the ideal spanned by the images of the
elements form $\R(\g)$. The subalgebra $\g\subset\so(p,q)$ is called
a Berger subalgebra if $L(\R(\g))=\g$. From the Ambrose-Singer
Theorem it follows that the holonomy algebra $\g\subset\so(p,q)$ of
a pseudo-Riemannian manifold of signature $(p,q)$ is a Berger
subalgebra.

A subalgebra $\g\subset\so(p,q)$ is called {\emph weakly
irreducible} if it does not preserve any proper non-degenerate
subspace of the pseudo-Euclidean space $\Real^{p,q}$. By the Wu
Theorem \cite{Wu}, any pseudo-Riemannian manifold whose holonomy
algebra is not weakly irreducible, can be decomposed (at least
locally) in the product of a flat pseudo-Riemannian manifold and of
pseudo-Riemannian manifolds with weakly irreducible holonomy
algebras. In particular, a pseudo-Riemannian manifold is  locally
indecomposable if and only if its holonomy algebra is weakly
irreducible. The following lemma is the well-known algebraic
counterpart of the Wu theorem.

\begin{lem}\label{lem Wu g} If a Berger subalgebra $\g\subset\so(p,q)$ is not weakly irreducible, then there
exists an orthogonal decomposition
$$\Real^{p,q}=V_0\oplus V_1\oplus\cdots\oplus V_r$$ into a direct sum of
pseudo-Euclidean subspaces and a decomposition
$$\g=\g_1\oplus\cdots\oplus \g_r$$ into a direct sum of ideals such that $\g_i$
annihilates $V_j$ if $i\neq j$ and $\g_i\subset \so(V_i)$ is a
weakly irreducible Berger subalgebra. \end{lem}

Consider the vector space
$$\R^\nabla(\g)= \left\{S\in (\Real^{p,q})^* \otimes\R(\g)\left|\begin{matrix}S_X(Y,Z)+S_Y(Z,X)+S_Z(X,Y)=0\\
\text{for all  } X,Y,Z\in \Real^{p,q}
\end{matrix}\right\}\right..$$ If a Berger algebra $\g\subset\so(p,q)$ satisfies $\R^\nabla(\g)=0$,
then $\g$ is called {\emph a symmetric Berger algebra}. Any
pseudo-Riemannian manifold with such holonomy algebra  is
automatically locally symmetric.
 The list of irreducible  holonomy algebras of locally symmetric
pseudo-Riemannian manifolds can be found in \cite{Ber57}.

M.~Berger classified irreducible subalgebras $\g\subset\so(p,q)$
that satisfy $L(\R(\g))=\g$ and $\R^\nabla(\g)\neq 0$ in
\cite{Ber}. Later this classification was corrected \cite{Al1},
and it was proved that all the obtained algebras can be realized
as the holonomy algebras of pseudo-Riemannian manifolds
\cite{Bryant2,Joyce07}. Here is the list of these algebras:

\begin{center}
$\so(p,q)$, $\so(p,\Co)\subset\so(p,p)$,\\
$\u(r,s)\subset\so(2r,2s)$, $\su(r,s)\subset\so(2r,2s)$,\\
$\sp(r,s)\subset\so(4r,4s)$,
$\sp(r,s)\oplus\sp(1)\subset\so(4r,4s)$,\\
 $\sp(r,\Real)\oplus\sl(2,\Real)\subset\so(2r,2r)$,
 $\sp(r,\Co)\oplus\sl(2,\Co)\subset\so(4r,4r)$,\\
  $\spin(7)\subset\so(8)$, $\spin(4,3)\subset\so(4,4)$,
  $\spin(7,\Co)\subset\so(8,8)$,\\
$G_2\subset\so(7)$, $G^*_{2(2)}\subset\so(4,3)$,
$G^\Co_2\subset\so(7,7)$.
\end{center}

Consider now the Einstein condition. For a subalgebra
$\g\subset\so(p,q)$, consider the set
$$\R_1(\g)=\{R\in\R(\g)|\Ric(R)(X,Y)=g(X,Y),\,\,X,Y\in\Real^{p,q}\}.$$
This is an affine space with the corresponding vector space
$$\R_0(\g)=\{R\in\R(\g)|\Ric(R)=0\}.$$
Let $L(\R_1(\g))\subset\g$ be the vector subspace spanned by the
images of the elements of $R\in\R_1(\g)$. We say that a subalgebra
$\g\subset\so(p,q)$ is {\emph a Berger subalgebra of Einstein type}
if $L(\R_1(\g))=\g$. From the Ambrose-Singer Theorem it follows that
{\emph the holonomy algebra $\g\subset\so(p,q)$ of an Einstein not
Ricci-flat pseudo-Riemannian manifold of signature $(p,q)$ is a
Berger subalgebra of Einstein type.} Note that in the Einstein not
Ricci-flat case, in Lemma \ref{lem Wu g} it holds $V_0=0$ (see the
proof of Lemma \ref{lem vid h} below).

In Section \ref{sec B irr} we will explain the classification of
irreducible holonomy algebras of Einstein not Ricci-flat
pseudo-Riemannian manifolds.

Now, the problem is to classify weakly irreducible and not
irreducible holonomy algebras of Einstein not Ricci-flat
pseudo-Riemannian manifolds. Until now this has been done only for
the Lorentzian signature in \cite{GalEin}.

Let $\g\subset\so(p,q)$ be a weakly irreducible not irreducible
subalgebra. Then $\g$ preserves a totally isotropic subspace
$V\subset\Real^{p,q}$ of dimension $m$. Let $(p,q)=(m+r,m+s)$. Let
$p_1,...,p_m$ be a basis of $V$. Choose linearly independent vectors
$q_1,...,q_m$ such that $g(q_i,q_j)=0$ and $g(p_i,q_j)=\delta_{ij}$.
Note that there is no canonical choice of these vectors.  The vector
space spanned by $q_1,...,q_m$ can be identified with $V^*$. We will
denote $V$ by $\Real^m$. Let $L\subset\Real^{m+r,m+s}$ be the
orthogonal complement to $\Real^m\oplus\Real^{m*}$. This subspace
can be identified with $\Real^{r,s}$. Let $e_1,...,e_{r+s}$ be an
orthonormal basis of $\Real^{r,s}$. The basis
$p_1,...,p_m,e_1,...,e_{r+s},q_1,...,q_m$ is called a Witt basis of
$\Real^{m+r,m+s}$.
 The maximal subalgebra
$\so(m+r,m+s)_{\Real^m}\subset\so(m+r,m+s)$ preserving the
subspace $\Real^m\subset\Real^{m+r,m+s}$ has the following matrix
form:
$$\so(m+r,m+s)_{\Real^m}=\left.\left\{\left(\begin{array}{ccc}
B&-X^tE_{r,s}&C\\
0&A&X\\
0&0&-B^t\end{array}\right)\right|\begin{matrix}
B\in\gl(m,\Real),\,A\in\so(r,s),\\
\,X\in\Real^m\otimes\Real^{r,s},\,C\in\wedge^2\Real^m\end{matrix}\right\}.$$
Denote the element given by the above matrix by $(B,A,X,C)$. The
non-zero Lie brackets are the following:
$$[(B_1,A_1,0,0),(B_2,A_2,0,0)]=([B_1,B_2]_{\gl(m,\Real)},[A_1,A_2]_{\so(r,s)},0,0),$$
$$[(B,A,0,0),(0,0,X,C)]=(0,0,AX+XB^t,BC+CB^t),$$
$$[(0,0,X,0),(0,0,Y,0)]=(0,0,0,-X^tE_{r,s}Y+Y^tE_{r,s}X).$$
We see that $\gl(m,\Real)$ and $\so(r,s)$ are subalgebras in
$\so(m+r,m+s)_{\Real^m}$. The vector subspace
$$\N=\{(0,0,X,0)|X\in\Real^m\otimes\Real^{r,s}\}\subset\so(m+r,m+s)_{\Real^m}$$
is isomorphic to the $\gl(m,\Real)\oplus\so(r,s)$-module
$\Real^m\otimes\Real^{r,s}$. The ideal
$$\C=\{0,0,0,C|C\in\wedge^2\Real^m\}\subset\so(m+r,m+s)_{\Real^m}$$
is isomorphic to the $\gl(m,\Real)$-module $\wedge^2\Real^m$. Note
that $\N\zr\C\subset\so(m+r,m+s)_{\Real^m}$ is a solvable ideal.
We get the decomposition
$$\so(m+r,m+s)_{\Real^m}=\gl(m,\Real)\oplus\so(r,s)\zr(\N\zr\C).$$
We will consider projections with respect to this decomposition.

Any  subalgebra $\g\subset\so(m+r,m+s)$ preserving a totally
isotropic subspace of dimension $m$ is conjugated to a subalgebra
of $\so(m+r,m+s)_{\Real^m}$.

If the holonomy algebra of a pseudo-Riemannian manifold $(M,g)$ of
signature $(m+r,m+s)$ preserves a totally isotropic subspace of the
tangent space of dimension $m$, then locally on $M$ there exists the
Walker \cite{Walker} coordinates $v_1,...,v_m,x_1,...,x_n,
u_1,...,u_m$ ($n=r+s$) such that the metric $g$ is of the form
$$g=\sum_{a=1}^m2dv_adu_a+h+\sum_{a=1}^m\sum_{b=1}^n2A_{ab}dx_adu_b+\sum_{a,b=1}^m
H_{ab}du_adu_b,$$ where $h=\sum_{a,b=1}^nh_{ab}(x_1,...,x_n,
u_1,...,u_m)dx_adx_b$ is a family of pseudo-Riemannian metrics of
signature $(r,s)$ depending on the parameters $u_1,...,u_m$;
$A_{ab}$ are functions of $x_1,...,x_n, u_1,...,u_m$, and $H_{ab}$
are functions of all coordinates.

\section{The case of irreducible $\g\subset\so(p,q)$}\label{sec B irr}

The list of irreducible holonomy algebras $\g\subset\so(p,q)$ of
Einstein not Ricci-flat pseudo-Riemannian manifolds of signature
$(p,q)$ that are not locally symmetric is the following (see e.g.
\cite{Bryant2}): \begin{center}$\so(p,q)$,
$\so(p,\Co)\subset\so(p,p)$, $\u(r,s)\subset\so(2r,2s)$,\\
$\sp(r,s)\oplus\sp(1)\subset\so(4r,4s)$,
$\sp(r,\Real)\oplus\sl(2,\Real)\subset\so(2r,2r)$,\\
$\sp(r,\Co)\oplus\sl(2,\Co)\subset\so(4r,4r)$.\end{center}

Consider a simply connected pseudo-Riemannian symmetric space
$(M,g)$ with an irreducible holonomy algebra $\g\subset\so(p,q)$.
Let $\k$ be the Lie algebra of the Lie group of transvections of
$(M,g)$. The irreducibility of $\g$ implies that $\k$ is simple.
Note that if $\g\subset \so(n,n)$  preserves two complementary
isotropic subspaces and acts irreducibly on these subspaces, then
$\k$ is simple also in that case (we will consider this situation in
the next section).

Consider the symmetric decomposition
$$\k=\g\oplus\mathfrak{m},$$
where the subspace $\mathfrak{m}$ is identified with the tangent
space of $M$ at some point. The Ricci tensor at that point is given
by the restriction of the Killing form of $\k$ to $\mathfrak{m}$.
Consequently the Ricci tensor is non-degenerate. From that and the
Schur lemma it follows that the space of symmetric bilinear tensor
fields on $(M,g)$ is either one-dimensional or it is
two-dimensional. Since the Ricci tensor is parallel and
non-degenerate, in the first case the manifold is Einstein and not
Ricci-flat. In the second case, there is a parallel $g$-symmetric
complex structure $I$ on $(M,g)$. This implies that $q=p$ and
$\g\subset\so(p,\Co)\subset\so(p,p)$. Any parallel symmetric
bilinear tensor field is of the form
$$ag(\cdot,\cdot)+bg(I\cdot,\cdot).$$
If $a^2+b^2\neq 0$, then such a form is non-degenerate and it
defines the same Levi-Civita connection as $g$. Thus the metric on
$M$ may be chosen to be Einstein and not Ricci-flat.

\section{The case of $\g\subset\so(n,n)$ preserving two complementary isotropic
subspaces}\label{sec nn}

Let $\g\subset\so(n,n)$ and suppose that $\g$ preserves two
complementary isotropic subspaces $V,V'\subset\so(n,n)$. Using the
metric on $\Real^{n,n}$, the space $V'$ can be identified with the
dual space $V^*$. If we fix a Witt basis $p_1,...,p_n,
q_1,...,q_n$ such that $p_1,...,p_n\in V$, and $q_1,...,q_n\in
V^*$, then $\g$ is contained in the maximal subalgebra preserving
$V$ and $V^*$:
$$\g\subset\so(n,n)_{\Real^n}=\left.\left\{\left(\begin{array}{cc}
A&0\\0&-A^t\end{array}\right)\right|A\in\gl(n,\Real)\right\}\subset\so(n,n).$$
The matrices $A$ define a subalgebra of $\gl(n,\Real)$, which is
isomorphic to $\g$ and we write $\g\subset\gl(n,\Real)$.

Conversely, if we start with a subalgebra $\g\subset \gl(n,\Real)$
we may consider it as the subalgebra
$\g\subset\so(n,n)=\so(\Real^n\oplus\Real^{n*})$ that preserves
$\Real^n,\Real^{n*}\subset \Real^n\oplus\Real^{n*}$. The
pseudo-Euclidean metric on $\Real^n\oplus\Real^{n*}$ is given by the
pairing.

Holonomy algebras $\g\subset\so(n,n)$ of this type are studied in
\cite{BBnn,BI97}. Let us formulate some facts from these papers.
Let $R\in \R(\g\subset\so(n,n))$. If $X,Y\in\Real^n$, or
$X,Y\in\Real^{n*}$, then $R(X,Y)=0$. Next, if $Y\in \Real^{n*}$,
then the map $$X\in\Real^n\mapsto R(X,Y)$$ belongs to the first
prolongation $(\g\subset\gl(n,\Real))^{(1)}$ of the representation
$\g\subset\gl(n,\Real)$. Hence, if $\g\subset\so(n,n)$ is a Berger
subalgebra, then $(\g\subset\gl(n,\Real))^{(1)}\neq 0$. Similarly,
the elements of $ \R^\nabla(\g\subset\so(n,n))$ can be described
using the second prolongation $(\g\subset\gl(n,\Real))^{(2)}$.
Hence, if $\g\subset\so(n,n)$ is a Berger subalgebra, and
$(\g\subset\gl(n,\Real))^{(2)}=0$, then it is a symmetric Berger
subalgebra.

 Next, suppose that $\g$ is a reductive Lie algebra.
 Lemma \ref{lem Wu g} applied to this case can be formulated in
 the following way.

\begin{lem}\label{lem Wu g nn} Let $\g$ be a reductive Lie algebra and suppose
that $\g\subset\so(n,n)$ is a Berger subalgebra preserving two
complementary totally isotropic subspaces of $\Real^{n,n}$. If
$\g\subset\gl(n,\Real)$ is not irreducible, then there exists a
decomposition
$$\Real^{n}=V_0\oplus V_1\oplus\cdots\oplus V_r$$  and a decomposition
$$\g=\g_1\oplus\cdots\oplus \g_r$$ into a direct sum of ideals such that $\g_i$
annihilates $V_j$ if $i\neq j$, $\g_i\subset \gl(V_i)$ is
irreducible, and $\g_i\subset\so(V_i\oplus V_i^*)$ is a  Berger
subalgebra.
\end{lem}

Let us restrict the attention to the case when
$\g\subset\gl(n,\Real)$ is irreducible. In \cite{BBnn} it is shown
that if $\g\subset\so(n,n)$ is a Berger subalgebra, then
$\g\subset\gl(n,\Real)$ is either one of the following:
$$\gl(n,\Real),\quad \sl(n,\Real),\quad \sp(2m,\Real),\quad  \gl(m,\Co),\quad \sl(m,\Co),\quad \sp(2k,\Co),\quad
2m=4k=n,$$ or $\g\subset\so(n,n)$ is a  symmetric Berger subalgebra,
and it can be found in the list from \cite{Ber57}. This result may
be deduced from the list of irreducible subalgebras
$\g\subset\gl(n,\Real)$ with $\g^{(1)}\neq 0$ given in Table~B from
\cite{Bryant2}.

Let us apply now the Einstein condition. Let
$R\in\R(\g\subset\so(n,n))$ be as above; let $X\in V$, $Y\in V^*$,
then it is easy to calculate
\begin{equation}\label{Rictr}\Ric(R)(X,Y)=\tr(
\pr_{\gl(n,\Real)}R(X,Y)).\end{equation} This shows that if
$\R_1(\g)\neq 0$, then $\g$ is not contained in $\sl(n,\Real)$.
Consequently, if $\g\subset\gl(n,\Real)$ is irreducible, then $\g$
is either one of
$$\gl(n,\Real),\quad \gl(m,\Co)\,\,\, (n=2m),$$ or $\g\subset\so(n,n)$ is a
weakly irreducible symmetric Berger algebra. Note that, as in the
previous section, it can be shown that the corresponding symmetric
spaces may  be Einstein and not Ricci-flat. Moreover,
$\gl(n,\Real)\subset\so(n,n)$ and $\gl(m,\Co)\subset\so(2m,2m)$ may
appear as the holonomy algebras of symmetric spaces. Thus, we have
found all Berger algebras of Einstein type $\g\subset\so(n,n)$ such
that $\g\subset\gl(n,\Real)$ is irreducible. Moreover, all these
Berger algebras are holonomy algebras of Einstein not Ricci-flat
pseudo-Riemannian manifolds of signature $(n,n)$. We write down that
list in Table~1.

\begin{tab} Irreducible subalgebras $\g\subset\gl(n,\Real)=\gl(V)$ defining  Berger algebras of Einstein type
$\g\subset\so(n,n)$

\begin{tabular}{| c | c | c |}\hline
$\g$ & $V$ & \text{ restriction }\\ \hline
$\gl(m,\Real)$&$\Real^m$&$m\geq 1$\\
$\gl(m,\Co)$&$\Co^m$&$m\geq 1$\\ \hline
$\Real\oplus\so(p,q)$&$\Real^{p,q}$&$p+q\geq 3$\\
$\Co\oplus\so(p,\Co)$&$\Co^{p}$&$p\geq 3$\\ \hline
$\Real\oplus\sl(p,\Real)\oplus\sl(q,\Real)$&$\Real^{pq}$&$p\geq q\geq 2$, $(p,q)\neq (2,2)$\\
$\Co\oplus\sl(p,\Co)\oplus\sl(q,\Co)$&$\Co^{pq}$&$p\geq q\geq 2$, $(p,q)\neq (2,2)$\\
$\Real\oplus\sl(p,\mathbb{H})\oplus\sl(q,\mathbb{H})$&$\Real^{4pq}$&$p\geq
q\geq 1$, $(p,q)\neq (1,1)$\\ \hline
$\Real\oplus\sl(p,\Co)$&$\Real^{p^2}\simeq H_p(\Co)$&$p\geq 3$\\
\hline
$\gl(p,\Real)$&$\Real^{p(p+1)/2}\simeq S_p(\Real)$&$p\geq 3$\\
$\gl(p,\Co)$&$\Co^{p(p+1)/2}\simeq S_p(\Co)$&$p\geq 3$\\
$\gl(p,\mathbb{H})$&$\Real^{p(2p+1)}\simeq S_p(\mathbb{H})$&$p\geq
2$\\\hline
$\gl(p,\Real)$&$\Real^{p(p-1)/2}\simeq A_p(\Real)$&$p\geq 5$\\
$\gl(p,\Co)$&$\Co^{p(p-1)/2}\simeq A_p(\Co)$&$p\geq 5$\\
$\gl(p,\mathbb{H})$&$\Real^{p(2p-1)}\simeq A_p(\mathbb{H})$&$p\geq
3$\\ \hline $\Real\oplus\spin(5,5)$&$\Real^{16}$&\\
$\Real\oplus\spin(1,9)$&$\Real^{16}$&\\
$\Co\oplus\spin(10,\Co)$&$\Co^{16}$&\\\hline $\Real\oplus
E^1_6$&$\Real^{27}$&\\
$\Real\oplus
E^4_6$&$\Real^{27}$&\\
$\Co\oplus E^\Co_6$&$\Co^{27}$&\\ \hline
\end{tabular}
\end{tab}

Ending this section we note that from \eqref{Rictr} follows the
following general statement.

\begin{prop} Let $(M,g)$ be a pseudo-Riemannian manifold of signature
$(n,n)$ such that its holonomy algebra $\g\subset\so(n,n)$ preserves
two complementary totally isotropic subspaces of the tangent space.
Then $(M,g)$ is Ricci-flat if and only if
$\g\subset\sl(n,\Real)$.\end{prop}

\section{Some lemmas}

\begin{lem}\label{lem-zgirr} Let $\g\subset\so(r,s)$ be an irreducible Berger subalgebra of Einstein type, then $\g$ contains its centralizer in $\so(r,s)$. \end{lem}

{\emph Proof.}  This immediately  follows from the classification
of irreducible holonomy algebras $\g\subset\so(r,s)$ of Einstein
pseudo-Riemannian manifolds. $\Box$

\begin{lem}\label{lem-zgnn} Let $\g\subset\gl(n,\Real)$ be an irreducible subalgebra such that
$\g\subset\so(n,n)$ is a Berger subalgebra of Einstein type, then
$\g$ contains its centralizer in $\so(n,n)$. \end{lem}

{\emph Proof.} As the $\g$-module, the Lie algebra $\so(n,n)$
decomposes as $$\so(n,n)=\gl(n,\Real)\oplus
\wedge^2\Real^n\oplus\wedge^2\Real^{n*}.$$ The classification shows
that the centralizer of $\g$ in $\gl(n,\Real)$ is contained in $\g$.
Next, $\g$ contains $\id_{\Real^n}$ that acts as multiplication by 2
and -2 in $\wedge^2\Real^{n}$ and $\wedge^2\Real^{n*}$,
respectively. $\Box$

\begin{lem}\label{lem-Rgirr} Let $\g\subset\so(r,s)$ be an irreducible  subalgebra and
$\R_1(\g)\neq 0$, then
$\g\subset\so(r,s)$ is a Berger subalgebra of Einstein type. \end{lem}

{\emph Proof.} Since the representation $\g\subset\so(r,s)$ is
irreducible, the Lie algebra $\g$ is reductive. Consider the ideal
$\k=L(\R_1(\g))\subset\g$. Suppose that $\k\neq\g$. Let
$\tilde\k\subset\g$ be the complementary ideal. The subalgebra
$\k\subset\so(r,s)$ is a Berger subalgebra, and
$\R_1(\k)=\R_1(\g)\neq 0$. By Lemma \ref{lem-zgirr}, the ideal
$\tilde\k\subset\g$ is not commutative. If the semi-simple part of
$\g$ is simple, then the statement of the lemma follows from Lemma
\ref{lem-zgirr}. Assume that the semi-simple part of $\g$ is not
simple.

Consider the case when the complexified representation
$\g\otimes\Co\subset\so(r+s,\Co)$ is not irreducible. Then
$\g\otimes\Co$ preserves two complementary isotropic subspaces
$V,V'\subset\Co^{r+s}$, and the induced representations are
irreducible. Let
$$R\in\R(\g\otimes\Co\subset\so(r+s,\Co))=\R(\g\subset\so(r,s))\otimes\Co.$$
As in Section \ref{sec nn}, for each $Y\in V'$ it holds
$$R(Y,\cdot)\in (\g\subset\gl(V))^{(1)}.$$ This shows that
$$(\g\subset\gl(V))^{(1)}\neq 0.$$ Using the list of irreducible
representations of complex Lie algebras with non-trivial first
prolongations (see e.g. \cite{Bryant2}) and the assumption that the
 semi-simple part of $\g$ is not simple, we get that
$$\g\otimes\Co=\gl(n,\Co)\oplus\gl(m,\Co),\quad V=\Co^n\otimes\Co^m$$
for some $n$ and $m$. This implies that
$\g=\u(n_1,n_2)\oplus\u(m_1,m_2)$, which is a Berger algebra of
Einstein type.

Finally suppose that the complexified representation
$\g\otimes\Co\subset\so(r+s,\Co)$ is irreducible. The representation
of $\g\otimes\Co$ in $\Co^{r+s}$ is the tensor products of
irreducible representations of $\k\otimes\Co$ and
$\tilde\k\otimes\Co$. The subalgebra
$\k\otimes\Co\subset\so(r+s,\Co)$ is a Berger subalgebra that
preserves at least two non-degenerate subspaces of $\Co^{r+s}$ and
its representation in each of the invariant subspaces is faithful.
This is impossible by the complex version of Lemma \ref{lem Wu g}.
This proves the current lemma. $\Box$

Note that in \cite{C-S} a  statement similar to Lemma
\ref{lem-Rgirr} was proven, stating that if $\R(\g)\neq 0$ for an
irreducible subalgebra $\g\subset\so(n)$, then either $\g$ is the
holonomy algebra of a Riemannian manifold or $\g=\sp(m)\oplus\u(1)$,
$n=4m$.

\begin{lem}\label{lem-Rginn} Let $\g\subset\gl(n,\Real)$ be an irreducible subalgebra such that
for  $\g\subset\so(n,n)$ it holds $\R_1(\g)\neq 0$, then
$\g\subset\so(n,n)$ is a Berger subalgebra of Einstein type.
\end{lem}

{\emph Proof.} In Section \ref{sec nn} we have seen that the
condition $\R_1(\g\subset\so(n,n))\neq 0$ implies
$$(\g\subset\gl(n,\Real))^{(1)}\neq 0,\text{ and }
\g\not\subset\sl(n,\Real).$$ Irreducible subalgebras
$\g\subset\gl(n,\Real)$ with non-trivial first prolongation may be
found in \cite{Bryant2}. This allows us to conclude that
$\g\subset\so(n,n)$ is a Berger subalgebra of Einstein type.
$\Box$

\begin{lem}\label{lem smena basisa}
Let $p_1,...,p_m,e_1,...,e_n,q_1,...,q_m$ be a Witt basis of
$\Real^{m+r,m+s}$ ($r+s=n$) and
$$\xi=(\id_{\Real^m},0,X,C)\in\gl(V)\oplus\so(r,s)\zr(\Real^m\otimes\Real^{r,s}\zr\wedge^2\Real^m)
=\so(m+r,m+s)_{\Real^m}$$ with respect to this basis. Then there
exists  another Witt basis of $\Real^{m+r,m+s}$ with the same
$p_1,...,p_m$, with respect to that
$$\xi=(\id_{\Real^m},0,0,0).$$
\end{lem}

{\emph Proof.} Let $p_1,...,p_m,e_1,...,e_n,q_1,...,q_m$ be a Witt
basis of $\Real^{m+r,m+s}$.  Consider the new basis
$$p'_i=p_i,\quad e'_a=e_a+\sum_{i=1}^mD_{ia}p_i,\quad
q'_i=q_i+X_i+\sum_{j=1}^mA_{ji} p_j,\quad X_i\in E.$$ Let $A=B+C$
be the decomposition of the matrix $A$ into the symmetric and
skew-symmetric parts. The condition that we get again a Witt basis
is equivalent to the equalities
$$B_{ij}=-\frac{1}{2}g(X_i,X_j),\quad D_{ia}=-g(X_i,e_a).$$
Let $\eta=(\id_{\Real^m},0,0,0)$ with respect to the first basis,
then with respect to the second basis,
$\eta=(\id_{\Real^m},0,-X,-2C))$, where $X_{ai}=-(X_i,e_a)$. It is
clear that $X$ and $C$ can be chosen in arbitrary way. This shows
that starting with $\xi=(\id_{\Real^m},0,X,C)\in\g$, we may choose
a new basis in such a way that $\xi=(\id_{\Real^m},0,0,0)\in\g$. $\Box$

\section{The general case}

Let $(M,g)$ be an Einstein not Ricci-flat pseudo-Riemannian manifold
with the holonomy algebra $\g\subset\so(p,q)$. We would like to find
all possible $\g$. As it was explained above, the Wu theorem allows
us to assume that $\g\subset\so(p,q)$ is weakly irreducible. The
list of irreducible $\g\subset\so(p,q)$ is known, see Section
\ref{sec B irr} above. Thus we may assume that $\g\subset\so(p,q)$
is weakly irreducible and it preserves a degenerate vector subspace
of $\Real^{p,q}$.

\begin{theorem}\label{th gener}
Let $\g\subset\so(r+m,s+m)$ be a weakly irreducible not irreducible
holonomy algebra of an Einstein not Ricci-flat pseudo-Riemannian
manifold of signature $(r+m,s+m)$. Let $m$ be the maximal dimension
of isotropic $\g$-invariant subspace.  Then with respect to a proper
basis of $\Real^{r+m,s+m}$, $\g\subset\so(r+m,s+m)_{\Real^m}$ has
the form
$$\g=\f\oplus\h\zr \left(\oplus_{i=1}^k\oplus_{\alpha=1}^t
\N_{i\alpha}\zr\oplus_{1\leq i\leq j\leq k}\C_{ij}\right),$$ where

\begin{itemize}
\item[$\bullet$] $\h\subset\so(r,s)$ is a reductive subalgebra that is the holonomy algebra of an
Einstein not Ricci-flat pseudo-Riemannian manifold of signature
$(r,s)$; there exist an orthogonal decomposition
$$\Real^{r,s}=L_1\oplus\cdots \oplus L_t,$$ and the corresponding
decomposition
$$\h=\h_1\oplus\cdots\oplus \h_t$$ such that $\h_\alpha\subset\so(L_\alpha)$
is an irreducible holonomy algebra of an  Einstein not Ricci-flat
pseudo-Riemannian manifold;

\item[$\bullet$] $\f\subset\gl(m,\Real)$ is a subalgebra that admits a
decomposition
$$\f=\f_0\zr \hat \f,$$ where $\f_0\subset\f$ is a reductive
subalgebra and $\hat\f\subset\f$ is a solvable ideal; the subalgebra
$\f_0\subset\gl(m,\Real)$ defines  the decompositions
$$\Real^{m}=V_1\oplus\cdots \oplus V_k,$$
$$\f_0=\f_1\oplus\cdots\oplus \f_k$$ such that
$\f_i\subset\gl(V_i)$ is irreducible and the corresponding weakly
irreducible subalgebra $\f_i\subset\so(V_i\oplus V_i^*)$ is the
holonomy algebra of an Einstein not Ricci-flat pseudo-Riemannian
manifold of neutral signature; farther more,
$$\hat\f=\oplus_{1\leq i<j\leq k}\f_{ij},$$ where $\f_{ij}\subset V^*_j\otimes
V_i$ is a $\f_i\oplus\f_j$-submodule; thus $\f\subset\gl(m,\Real)$
has the following structure:
$$\f=\left\{\left.\left(\begin{array}{ccccc}
A_1 & A_{12}& \cdots& A_{1\,k-1} &A_{1k}\\
0 &A_2&\cdots&  A_{2\,k-1} &A_{2k}\\
\vdots&\vdots&\vdots&\vdots&\vdots\\
0&0&\cdots&0&A_k\end{array} \right)\right| A_i\in\f_i,\,A_{ij}\in
\f_{ij}\right\};$$

\item[$\bullet$] $\N_{i\alpha}$ is an $\f_i\oplus \h_\alpha$-submodule of
$V_i\otimes L_\alpha$; for each $i$, there exists $\alpha$ such that
$\N_{i\alpha}\neq 0$;

\item[$\bullet$] $\C_{ij}$ is an $\f_i\oplus\f_j$-submodule of  $V_i\wedge V_j$
if $i\neq j$, and $C_{ii}\subset \wedge^2 V_i$ is an
$\f_i$-submodule;

\item[$\bullet$] it holds
\begin{align}\label{relations1}[\f_{li},\f_{ij}]&\subset\f_{lj},\\ \label{relations2}
[\f_{li},\C_{ij}]&\subset\C_{lj},\\ \label{relations3}
[\f_{li},\N_{i\alpha}]&\subset\N_{l\alpha},\\ \label{relations4}
[\N_{i\alpha},\N_{j \alpha}]&\subset\C_{ij}.\end{align}
\end{itemize}
\end{theorem}

{\bf Proof of Theorem \ref{th gener}.}

In the proof of the theorem we will use only the property
$\R_1(\g)\neq 0$.

 The following lemma gives
the form of the projection $\h=\pr_{\so(r,s)}\g$.

\begin{lem}\label{lem vid h} There exists an orthogonal decomposition $$\Real^{r,s}=L_1\oplus\cdots \oplus L_t,$$
and the corresponding decomposition $$\h=\h_1\oplus\cdots\oplus
\h_t$$ such that $\h_i\subset\so(L_i)$ is an irreducible holonomy
algebra of an Einstein not Ricci-flat pseudo-Riemannian manifold.
\end{lem}

{\emph Proof.} Since $\Real^m\subset\Real^{r+m,s+m}$ is the maximal
totally isotropic subspace preserved by $\g$, $\h$ does not preserve
any totally isotropic subspace of $\Real^{r,s}$, and consequently
$\h$ does not preserve any degenerate subspace of $\Real^{r,s}$.
Consequently, if $\h$ preserves a proper vector subspace of
$\Real^{r,s}$, then this subspace is non-degenerate, and $\g$
preserves also its orthogonal complement. This shows that $\h$ is a
reductive Lie algebra. We get an $\h$-invariant orthogonal
decomposition
$$\Real^{r,s}=U_0\oplus U_1\oplus\cdots\oplus U_k$$ such that $\h$
annihilates $U_0$ and  the induced representation of $\h$ on $U_i$
is irreducible for $1\leq i\leq k$. For each $i$, let
$\h_i\subset\h$ be the ideal annihilating $U_i^\bot$, i.e.
$\h_i\subset\so(U_i)$. We obtain the decomposition
$$\h=\h_1\oplus\cdots\oplus\h_k\oplus\bar \h,$$ where $\bar
\h\subset \h$ is the complementary ideal to $\h_1\oplus\cdots\oplus
\h_k$. The elements of $\bar \h$ act simultaneously in several
$U_i$. Note that $\h_i\subset\so(U_i)$ a priori must not be
irreducible.

For $R\in \R_1 (\g)$, let
$$\tilde R=\pr_{\so(r,s)}\circ R|_{\wedge^2\Real^{r,s}}.$$
Applying the Bianchi identity for the tensor $R$ to vectors from
$\Real^{r,s}$ and taking the projection to $\Real^{r,s}$, we see
that $\tilde R\in\R (\h)$.  Let $\tilde \h\subset\h$ be the
subalgebra generated by the images of such tensors $\tilde R$.
Clearly, $\tilde \h\subset\so(r,s)$ is a Berger subalgebra.

Let again $R\in \R_1 (\g)$ and let $X,Y\in\Real^{r,s}$. It holds
$$\Ric(R)(X,Y)=\sum_{a=1}^m g(R(p_a,X)Y,q_a)+\sum_{a=1}^{r+s}
g(R(e_a,X)Y,e_a)g(e_a,e_a)+\sum_{a=1}^m g(R(q_a,X)Y,p_a).$$ Since
$R$ takes values in $\g$, and $\g$ preserves the spaces $\Real^m$
and $\Real^m\oplus\Real^{p,q}$, we get $g(R(q_a,X)Y,p_a)=0$. Using
the standard property of the curvature tensor, we get
$$g(R(p_a,X)Y,q_a)=g(R(Y,q_a)p_a,X)=0.$$
This shows that
$$g(X,Y)=\Ric(R)(X,Y) =\Ric(\tilde R)(X,Y),$$ and $\tilde R\in\R_1
(\h)$.
We claim that $\tilde \h$ does not annihilate any non-degenerate
subspace of $\Real^{r,s}$. Indeed, if $\tilde\h$ annihilates a
non-zero vector $Z\in \Real^{r,s}$, then for any $X\in \Real^{r,s}$
it holds
$$g(X,Z)=\Ric(\tilde R)(X,Z)=\sum_{a=1}^{r+s}
g(R(e_a,X)Z,e_a)g(e_a,e_a)=0,$$ which is impossible.

Since $\tilde \h\subset\so(r,s)$ is a Berger subalgebra, we obtain
the decompositions
$$\Real^{r,s}=W_1\oplus\cdots \oplus W_k,$$
 $$\tilde \h=\tilde\h_1\oplus\cdots\oplus\tilde \h_k,$$
where each $\tilde\h_i\subset\so(W_i)$ is a weakly irreducible
Berger subalgebra.

Since $\tilde \h\subset \h$, each $U_i$ is the direct sum  of some number of $W_j$,
and the corresponding $\tilde \h_j$ are contained in $\h_i$. In particular, all $\h_i$ are non-trivial and do not annihilate any non-trivial subspace of $U_i$, and $U_0=0$.
More over, $\R_1(\h_i)\neq 0$.
 Consider the ideal $\t_i=L(\R(\h_i))\subset\h_i$.
Let $\bar\h_i\subset\h_i$ be the complementary ideal, then  $$\h_i=\t_i\oplus \bar\h_i.$$
Clearly,
$\t_i\subset\so(U_i)$ is a Berger subalgebra.

Since $U_i$ is the direct sum of some number of $W_j$, $\t_i$
contains the corresponding $\tilde\h_j$, and
$\tilde\h_j\subset\so(W_j)$ is weakly irreducible, we see that
$\t_i\subset\so(U_i)$  does not annihilate any non-trivial subspace
of $U_i$.
 We get the decompositions
$$U_i=U_{i1}\oplus\cdots\oplus U_{it_1},\quad \t_i=\t_{i1}\oplus\cdots\oplus\t_{it_i},$$
where $\t_{i\alpha}\subset\so(U_{i\alpha})$ is a weakly irreducible
Berger subalgebra of Einstein type.

By the construction, each $\t_{i\alpha}$ is a reductive Lie algebra,
hence either $\t_{i\alpha}\subset\so(U_{i\alpha})$ is irreducible,
or it preserves two  complementary totally isotropic subspaces. By
Lemmas \ref{lem-zgirr} and \ref{lem-zgnn}, the induced action of
$\bar \h_j$ and $\bar\h$ on each $U_{i\alpha}$ is trivial, i.e.
$\bar\h=0$ and all $\bar\h_j=0$. Hence,
$\h=\oplus_{i,\alpha}\t_{i\alpha}$. Since $\h$ does not preserve any
degenerate subspace of $\Real^{r,s}$, each
$\t_{i\alpha}\subset\so(U_{i\alpha})$ is irreducible. We obtain the
required decompositions $\Real^{r,s}=\oplus_{i,\alpha}U_{i\alpha}$,
$\h=\oplus_{i,\alpha}\t_{i\alpha}$. $\Box$

Now we find the form  of the projection $\f=\pr_{\gl(m,\Real)}\g$.

\begin{lem}\label{lem vid f}
 There exists a decomposition $\f=\f_0\zr \hat \f$, where $\f_0\subset\f$ is a reductive
subalgebra and $\hat\f\subset\f$ is a solvable ideal. It holds
$\f_0\subset\g$. Next, there is a decomposition
\begin{equation}\label{dec Rm}\Real^m=V=V_1\oplus\cdots \oplus V_k,\end{equation} and
the corresponding decomposition
\begin{equation}\label{dec f0}\f_0=\f_1\oplus\cdots\oplus \f_k\end{equation} such that
$\f_i\subset\gl(V_i)$ is irreducible and the corresponding weakly
irreducible subalgebra $\f_i\subset\so(V_i\oplus V_i^*)$ is the
holonomy algebra of an Einstein not Ricci-flat pseudo-Riemannian
manifold of neutral signature.

Moreover, as the vector space, $$\hat\f=\oplus_{1\leq i<j\leq
k}\f_{ij},$$ each $\f_{ij}$ takes $V_j$ to $V_i$ and annihilates
$V_l$, $l\neq j$.
\end{lem}

{\emph Proof.} Suppose that the representation of $\f$ on $V$
preserves a vector subspace $V_1\subset V$. We may assume that $V_1$
does not contain any proper invariant subspace, i.e. the induced
representation of $\f$ on $V_1$ is irreducible. Let $V'_1\subset V$
be any complementary subspace, i.e. $V=V_1\oplus V'_1$. The matrices
of the elements from $\f$ are of the form
$$\left(\begin{array}{cc} A_1&A_{12}\\0&A_2\end{array}\right).$$
That is, $$\f\subset\gl(V_1)\oplus\gl(V'_1)\zr (V'_1)^*\otimes
V_1.$$ Let $\f_1\subset\gl(V_1)$ be the projection of $\f$ with
respect to that decomposition. By the construction,
$\f_1\subset\gl(V_1)$ is irreducible. Consequently, $\f_1$ is a
reductive Lie algebra. Consider the decomposition
$$\Real^{m+r,m+s}=V_1\oplus L'\oplus V^*_1,$$ where  $L'=V'_1\oplus L\oplus
(V'_1)^*$. We get
$$\g\subset\gl(V_1)\oplus\so(L')\zr (V_1\otimes
L'\zr\wedge^2V_1).$$ Let $R\in \R_1 (\g)$ be a non-zero element.
Consider
$$\tilde R=\pr_{\gl(V_1)}\circ R|_{V_1\otimes V_1^*}.$$
As in the proof of Lemma \ref{lem vid h}, it is not hard to show
that $\tilde R\in\R_1(f_1\subset\so(V_1\oplus V_1^*))$. Also it is
easy to see that if $X\in V_1$, and $Y\in V_1^*$, then
$$g(X,Y)=\Ric(R)(X,Y)=\Ric(\tilde R)(X,Y).$$
This implies that $\tilde R\neq 0$. By Lemma \ref{lem-Rginn},
$\f_1\subset\gl(V_1)$ is a Berger subalgebra of Einstein type; in
Section \ref{sec nn} we noted that any such algebra is the holonomy
algebra of an Einstein  not Ricci-flat pseudo-Riemannian manifold.
From results of Section \ref{sec nn} it follows that $\f_1$ contains
$\id_{V_1}$.

We claim that $\pr_{\gl(V_1)\oplus\so(L')}\g=\f_1\oplus
\pr_{\so(L')}\g$. In other words, the ideal $$\k=\f_1\cap
\pr_{\gl(V_1)\oplus\so(L')}\g\subset\f_1$$ coincides with $\f_1$.

Using the Bianchi identity it is easy to get
\begin{equation}\label{R Rm Rrs} \pr_{\gl(V_1)}\circ
R|_{\wedge^2 L'}=0, \quad \pr_{\so(L')}\circ R|_{V_1\otimes
V_1^*}=0.\end{equation} The first of these equalities shows that the
above defined tensor $\tilde R$ takes values in $\k$. In particular,
$\R_1(\k\subset\gl(V_1\oplus V_1^*))\neq 0$.
Let $\tilde \k$ be the ideal complementary to $\k$ in $\f_1$.

Suppose that $\k\subset\gl(V_1)$ is irreducible. Since
$\R_1(\k\subset\gl(V_1\oplus V_1^*))\neq 0$, we see that by Lemma
\ref{lem-Rginn}, $\k\subset\gl(V_1\oplus V_1^*)$ is a Berger
subalgebra of Einstein type. Note that $\tilde
\k\subset\z_{\gl(V_1)}\k$. By Lemma~\ref{lem-zgnn}, $\tilde
\k\subset \k$. Thus, $\tilde\k=0$.

Suppose that $\k\subset\gl(V_1)$ is not irreducible. Consider the
ideal $$\t=L(\R(\k\subset\so(V_1\oplus V_1^*)))\subset\k.$$ We see
that $\t$ is not trivial, $\t\subset\so(V_1\oplus V_1^*)$ is a
Berger subalgebra, and $$\R_1(\t\subset\so(V_1\oplus V_1^*))\neq
0.$$ By the above arguments with the Ricci tensor, $\t$ does not
annihilate any non-trivial subspace of $V_1$. We obtain the
decompositions
$$V_1=U_1\oplus\cdots\oplus U_a, \quad \t=\t_1\oplus\cdots\oplus\t_a,$$
where for each $i$, $\t_i\subset\gl(U_i)$ is irreducible and
$\t_i\subset\so(U_i\oplus U_i^*)$ is a Berger subalgebra. By the
argument with the Ricci tensor, $\R_1(\t_i\subset\so(U_i\oplus
U_i^*))\neq 0$. By Lemma \ref{lem-Rginn}, $\t_i\subset\so(U_i\oplus
U_i^*)$ is a Berger subalgebra of Einstein type. This implies that
$\t_i$ contains $\id_{U_i}$. Then $\f_1$ contains the elements
$\id_{U_1}$ and $\id_{U_2}$, which by the construction commute with
$\f_1$. This contradicts the Schur Lemma applied to the irreducible
representation $\f_1\subset\gl(V_1)$. Thus, $\f_1=\k$ and the claim
is proved.

From Lemma \ref{lem-zgnn} it follows that $\id_{V_1}\in\f_1$. From
this and the above claim it follows that
 $(\id_{V_1},0,X,C)\in\g$ for some $X\in V_1\otimes L'$ and $C\in\wedge^2V_1$.
 Lemma \ref{lem smena basisa} shows that we may change the subspaces $L',V_1^*\subset\Real^{m+r,m+s}$
 in such a way that $(\id_{V_1},0,0,0)\in\g$. Note that under this change the projection $\pr{\gl(V_1)}\g=\f_1$
 does not change; the new space $L'$ is isomorphic to the old one, and the projection $\pr_{\so(L')}\g$
 does not change if we use this isomorphism as the identification.
 Let $A\in\f_1$. Then $(A,,0,X,C)\in\g$ for some $X\in V_1\otimes L'$ and $C\in\wedge^2V_1$. Taking the Lie brackets
  of this element with $(\id_{V_1},0,0,0)\in\g$, we get that $(0,0,X,2C)\in\g$. Taking the Lie brackets again,
 we get $(0,0,X,4C)\in\g$. We conclude that $\f_1\subset\g$. Similarly, $\pr_{\so(L')}\g\subset\g$, i.e.
 $\pr_{\so(L')}\g=\g\cap\so(L')$; and $$\g=\f_1\oplus(\g\cap\so(L'))\zr((\g\cap (V_1\otimes L'))\zr(\g\cap\wedge^2V_1)).$$

We consider now the intersection
$$\g\cap\so(L')\subset\gl(V')\oplus\so(L')\zr(V'\otimes
L'\zr\wedge^2 V_1).$$ As in the proof of Lemma \ref{lem vid h}, it
is easy to show that $\R_1(\g\cap\so(L'))\neq 0$. Consider the
projection $\pr_{\gl(V')} (\g\cap\so(L'))\subset\gl(V')$. If it is
not irreducible,  then it preserves a subspace $V_2\subset V'$
such that the induced representation on $V_2$ is irreducible. We
choose a complementary subspace $V_2'\subset V'$. Then we apply
the above consideration to this settings. We will get a similar
result for $\g\cap\so(L')$ as for $\g$ above. In particular,
$\f_2=\pr_{\gl(V_2)}(\g\cap\so(L'))\subset\gl(V_2)$ is irreducible and it
satisfies the same properties as $\f_1$.

We continue this process. At the last stage we will see  that $\f$
acts irreducibly on $V'_{k-1}$ for some $k$. We set
$V_k=V'_{k-1}$.  We get the $\f$-invariant decomposition
$$V=V_1\oplus\cdots\oplus V_k$$
and irreducible subalgebras $\f_i\subset\gl(V_i)$ that are
contained in $\f$ and $\g$; in particular, the  Lie algebras
$\f_i$ are reductive.

The representation of $\f$ on $\Real^m$ is given by matrices of
the form $$A=\left(\begin{array}{cccc}
A_1 & *& \cdots &*\\
0 &A_2&\cdots& *\\
\vdots&\vdots&\vdots&\vdots\\
0&0&\cdots&A_k\end{array} \right),\quad A_i\in\f_i\subset\g.$$ Let
$$\f_0=\f_1\oplus\cdots\oplus\f_k\subset\f.$$
Since the Lie algebra $\f_0$ is reductive, there exists a
complementary subspace $\hat\f\subset\f$ with respect to the adjoint
representation of $\f_0$ on $\f$. The above consideration shows that
the elements of $\hat\f$ are given by the stars in the matrix $A$.
In particular, $\hat\f\subset\f$ is an ideal. Clearly,
$\hat\f\subset \oplus_{1\leq i<j\leq k}V^*_j\otimes V_i$. We get the
semidirect sum $\f=\f_0\zr\f_1$. Consider the representation of
$\f_0$ on $\oplus_{1\leq i<j\leq k}V^*_j\otimes V_i$; on each
subspace $V_{ij}=V^*_j\otimes V_i$ it is given by the tensor product
of the representations $\f_j\subset\gl(V_j^*)$ and
$\f_i\subset\gl(V_i)$. Since all subalgebras $\f_i\subset\gl(V_i)$
are irreducible, we see that different $V_{ji}$ and $V_{ls}$ do not
contain isomorphic $\f_0$-modules. This shows that $\hat
\f=\oplus_{1\leq i<j\leq k}\f_{ij}$, where $\f_{ij}=\hat\f\cap
V_{ij}$. This proves the lemma. $\Box$

Now we consider $\pr_{\gl(m,\Real)\oplus\so(r,s)}\g$.

\begin{lem}\label{lem f+h} It holds
$\pr_{\gl(m,\Real)\oplus\so(r,s)}\g=\f\oplus\h$.\end{lem}

{\emph Proof.} Let $\k=\pr_{\gl(m,\Real)\oplus\so(r,s)}\g$. The
facts that $\f_0\subset\g$, $\h$ is a reductive Lie algebra, and
$\hat\f$ does not contain non-zero elements commuting with $\f_0$,
imply the decomposition
$$\k=\f_0\oplus \hat\f\oplus \h,$$
i.e. $\k=\f\oplus\h$.
 $\Box$

\begin{lem} It holds $\f\oplus\h\subset\g$. \end{lem}

{\emph Proof.} We have already seen that, for  each $i$, it holds
$\id_{V_i}\in\f_i\subset\g$. Hence, $\id_{\Real^m}\in\f_0\subset\g$.
Let $A\in \f$. From Lemma \ref{lem f+h} it follows that
$(A,0,X,C)\in\g$ for some $X\in\N$ and $C\in\C$. Taking the Lie
brackets of $\xi$ with that element, we get $(0,0,X,2C)\in\g$.
Taking the Lie brackets again, we get $(0,0,X,4C)\in\g$. We conclude
that $(A,0,0,0)\in\g$, i.e. $\f\subset\g$. Similarly, $\h\subset\g$.
$\Box$

From the last two lemmas it follows that
$\g=\f\oplus\h\zr(\g\cap(\N\zr\C)).$ As in the previous lemma, we
may show that if $(0,0,X,C)\in\g$, then $(0,0,X,0)\in\g$ and
$(0,0,0,C)\in\g$, that is, $$\g\cap(\N\zr\C)=(\g\cap \N)\zr
(\g\cap\C).$$

Next, as $\f_0\oplus\h$-module, $\N=\oplus_{i,\alpha}V_i\otimes
L_\alpha$.  It is clear that these $\f_0\oplus\h$-modules are
pairwise non-isomorphic, hence, $\g\cap\N=
\oplus_{i,\alpha}\N_{i\alpha}$, where  $\N_{i\alpha}=\g\cap (
V_i\otimes L_\alpha)$. Similarly, as the $\h$-module,
$\C=\oplus_{i\leq j}V_i\wedge V_j$, and $\g\cap\C= \oplus_{i\leq
j}\C_{ij}$, where  $\C_{ij}=\g\cap (V_i\wedge V_j)$.

Finally, if for some $\alpha$, all $\N_{i\alpha}$ are zero, then
$\g$ preserves the non-degenerate subspace
$L_\alpha\subset\Real^{r+m,s+m}$, this gives a contradiction,
since $\g$ is weakly irreducible. The theorem is proved. $\Box$

\vskip0.3cm

Let us now show that there are few possibilities for each $\f_{ij}$,
$\N_{i\alpha}$ and $\C_{ij}$. Recall that the tensor product of two
irreducible representations of Lie algebras in complex vector spaces
is always irreducible. It turns out that this is not the case for
the tensor product of real representations and we prove the
following (probably known) lemma.

\begin{lem}
Let $\h_i\subset\gl(V_i)$, $i=1,2$, be irreducible representations
of real Lie algebras in real vector spaces. If the representation
of $\h_1\oplus \h_2$ in $V_1\otimes V_2$ is not irreducible, then
 each vector space $V_i$ admits a complex structure $J_i$ commuting
with $\h_i$.
Moreover, $V_1\otimes V_2$ is the direct sum of two irreducible
subspaces of the form $V_1\otimes_\Co V_2$, where one considers
the tensor product   either  of the complex vector spaces
$(V_1,J_1)$ and $(V_2,J_2)$, or the tensor product of the complex
vector spaces $(V_1,-J_1)$ and $(V_2,J_2)$.
\end{lem}

{\emph Proof.} Suppose that the tensor product $V_1\otimes V_2$ is
not irreducible, and $V_1$ does not admit a complex structure
commuting with $\h_1$. Since the product $V_1\otimes V_2$ is not
irreducible, there exists a linear automorphism $K$ of $V_1\otimes
V_2$ commuting with $\h_1\oplus \h_2$ and such that $K^2=\id\neq
\pm K$. We may write $K$ in the form $$K=\sum_\alpha
A_\alpha\otimes B_\alpha,$$ where $A_\alpha\in\gl(V_1)$,
$B_\alpha\in\gl(V_2)$, and the vectors $\{B_\alpha\}$ are linearly
independent. Let $A\in \h_1$, $x\in V_1$ and $y\in V_2$, then
$$K(Ax\otimes y)=K(A(x\otimes y))=A(K(x\otimes y)),$$
i.e. $$\sum_\alpha A_\alpha Ax\otimes B_\alpha y=\sum_\alpha
AA_\alpha x\otimes B_\alpha y.$$ This implies that $\sum_\alpha
[A,A_\alpha]\otimes B_\alpha=0$ and $[A,A_\alpha]=0$. We conclude
that each $A_\alpha$ commutes with $\h_1$. Since there are no
complex structures on $V_1$ commuting with $\h_1$, we see that
from the Schur lemma it follows that $A_\alpha=c_\alpha
\id_{V_1}$, $c_\alpha\in\Real$. Thus,
$$K=\id_{V_1}\otimes B,\quad B\in\gl(V_2).$$
Since $K^2=\id_{V_1\otimes V_2}$, we get that $B^2=\id_{V_2}$.
Since $B$ commutes with $\h_2$, and $\h_2\subset\gl(V_2)$ is
irreducible, we see that $B=\pm \id_{V_2}$. Consequently,
$K=\pm\id_{V_1\otimes V_2}$ and we get a contradiction. Thus, each
$V_i$ admits a complex structure $J_i$ commuting with $\h_i$.

The two vector subspaces in $V_1\otimes V_2$ described in the
statement of the lemma are invariant, irreducible (as tensor
products of irreducible complex subspaces), their intersection is
trivial, and the sum of their real dimensions equals to the
dimension of $V_1\otimes V_2$, i.e. $V_1\otimes V_2$ is the direct
sum of these subspaces. $\Box$

We see now that there are only three possibilities for each space
$\f_{ij}$: it may be trivial; it may coincide with $V_j^*\otimes
V_i$; it may coincide with $V_j^*\otimes_\Co V_i$. In the last case,
there exist complex structures on $V_i$ and $V_j$ commuting with
$\f_i$ and $\f_j$, respectively. This means that $m_i$ and $m_j$ are
even, $$\f_i\subset\gl(m_i/2,\Co)\subset\gl(m_i,\Real),\quad
\f_j\subset\gl(m_j/2,\Co)\subset\gl(m_j,\Real),$$
 and $$\f_i\subset\u(m_i/2,m_i/2)\subset\so(m_i,m_i),\quad \f_j\subset\u(m_j/2,m_j/2)\subset\so(m_j,m_j)$$ are the holonomy algebras of
 Einstein not Ricci-flat
pseudo-K\"ahlerian manifolds.  Similarly, there are three
possibilities for each space $\N_{i\alpha}$ and $\C_{ij}$ with
$i\neq j$. In order to find the spaces $\C_{ii}$ one should, for
each possible $\f_i$, decompose the $\f_i$-module $\wedge^2V_i$ into
the direct sum of irreducible submodules. The space $\C_{ii}$ is
either trivial, or it is the direct sum of some of these submodules.

Let us explain now, how to list the Lie algebras $\g$ obtained in
Theorem \ref{th gener}. Fix a signature $(p,q)$, $p,q\geq 1$. Fix
numbers $r,s,m$ such that $p=r+m$, $q=s+m$, $m\geq 1$. Choose
decompositions $$m=m_1+\cdots+m_t,\quad t\geq 1,$$
$$r=r_1+\cdots+r_k,\quad s=s_1+\cdots+s_k,\quad k\geq 1.$$
Fix irreducible subalgebras $\h_\alpha\subset\so(r_\alpha,s_\alpha)$
and $\f_i\subset\gl(m_i,\Real)$ as in Sections \ref{sec B irr} and
\ref{sec nn}, respectively. Choose $\f_i\oplus \f_j$-modules
$\f_{ij}$ such that all conditions \eqref{relations1} are satisfied.
Choose $\f_i\oplus \h_\alpha$-modules $\N_{i\alpha}$ such that
\eqref{relations3} are satisfied. Choose $\f_i\oplus \f_j$-modules
$\C_{ij}$ and $\f_i$-modules $\C_{ii}$ that satisfy
\eqref{relations2} and \eqref{relations4}. We obtain a Lie
subalgebra $\g\subset\so(r+m,s+m)$. One may check if the obtained
subalgebra is weakly irreducible using the following proposition.

\begin{prop} A Lie subalgebra $\g\subset\so(r+m,s+m)$ described in
Theorem \ref{th gener} is not weakly irreducible if and only if
there are subsets $T_1,T_2\subset\{1,...,t\}$,
$K_1,K_2\subset\{1,...,k\}$ such that
$$T_1\cap T_2=\emptyset,\quad T_1\cup T_2=\{1,...,t\},\quad K_1\cap K_2=\emptyset,\quad
K_1\cup K_2=\{1,...,k\},$$ one of the subsets
$T_1\subset\{1,...,t\}$ and $K_1\subset\{1,...,k\}$ is proper, and
for each $i\in T_1$, $j\in T_2$, $\alpha\in K_1$, $\beta\in K_2$, it
holds $$\f_{ij}=0,\quad \N_{i\beta}=0,\quad \C_{ij}=0.$$
\end{prop}

{\bf Proof of the proposition.} If the above sets of indices exist,
then $\g$ preserves the proper non-degenerate vector subspace
$$\oplus_{i\in K_1} (V_i\oplus V^*_i) \bigoplus \oplus_{\alpha\in
T_1} L_\alpha\subset \Real^{r+m,s+m}.$$

If $\g\subset\so(r+m,s+m)$ is not weakly irreducible, then it can be
decomposed into the direct sum of two Berger algebras of Einstein
type. Applying Theorem \ref{th gener} to both algebras, we will get
the required sets of indices. $\Box$

We are left with the problem to construct an Einstein not Ricci-flat
pseudo-Riemannian manifold with each of the obtained weakly
irreducible holonomy algebras $\g\subset\so(r+m,s+m)$. For manifolds
of signatures $(1,n)$ and $(2,n)$ this is done in Sections \ref{sec
Lor} and \ref{sec m=2} below. The general case is discussed in
Section \ref{sec Rem}.

\section{Lorentzian manifolds}\label{sec Lor}  From Theorem \ref{th gener} it follows that any weakly
irreducible not irreducible holonomy algebra $\g\subset\so(1,n+1)$
of an Einstein not Ricci-flat Lorentzian manifold of dimension $n+2$
has the form
$$\g=\gl(1,\Real)\oplus\h\zr\Real^n,$$ where $\h\subset\so(n)$ is
the holonomy algebra of an Einstein not Ricci-flat Riemannian
manifold. This result for the first time was obtained in
\cite{GalEin}, for that was used the classification from
\cite{Leistner}. Note that the result of the present paper provides
an independent proof of the results from \cite{Leistner} in the
Einstein case. In \cite{GalEin}, the just described algebras $\g$
are realized as the holonomy algebras of Einstein not Ricci-flat
Lorentzian manifolds. Let $\Lambda\neq 0$. The required metric is
the following
$$g=2dvdu+h+(\Lambda v^2+H_0)(du)^2,$$ where $v,x^1,...,x^n,u$ are
coordinates on an open subset of $\Real^{n+2}$,
$$h=h_{ij}(x^1,...,x^n)dx^idx^j$$ is an Einstein Riemannian metric
with the holonomy algebra $\h\subset\so(n)$ and cosmological
constant $\Lambda$, and $H_0=H_0(x^1,...,x^n)$ is a function
satisfying $\Delta_h H_0=0$, where $\Delta_h$ is the Laplacian of
the metric $h$. It is required that $H_0$ has non-zero Hessian in
order make the metric to be indecomposable.

The idea of this construction is the following. Let $g_0$ be  the
metric $g$ with $H_0$ set to zero. Then $g_0$ is the direct product
of the Einstein Lorentzian manifold of dimension 2 with the holonomy
algebra $\gl(1,\Real)=\so(1,1)$ preserving two complement isotropic
lines and of a Riemannian manifold with the holonomy algebra $\h$;
the cosmological constant of each of these manifolds is $\Lambda$.
The function $H_0$ is used in order to curve the product metric and
to add the subspace $\Real^n$ to the holonomy algebra. The curvature
tensor of $g$ equals to the curvature tensor of $g_0$ with the
additional term that takes values in $\Real^n\subset\g$. This
additional term is given by a symmetric endomorphism of the
Euclidean space $\Real^n$, the Ricci tensor of the additional term
must be zero, and this gives the equation $\Delta_h H_0=0$.


\section{Pseudo-Riemannian manifolds of index 2}\label{sec m=2}

In this section we give a complete classification of the holonomy
algebras of Einstein not Ricci-flat pseudo-Riemannian manifolds of
signature $(2,n)$, $n\geq 2$.

We start with the signature $(2,2)$.

\begin{theorem}\label{th 22} Weakly irreducible not irreducible
holonomy algebras of Einstein not Ricci-flat pseudo-Riemannian
manifolds of signature $(2,2)$ are exhausted by the following
subalgebras of $\so(2,2)$:

{\bf 1)} $\g_1=(\gl(1,\Real)\oplus\gl(1,\Real))\zr \wedge^2\Real^2$.

{\bf 2)} $\g_2=\gl(2,\Real)$;

{\bf 3)} $\g_3=\gl(2,\Real)\zr\wedge^2\Real^2$;

{\bf 4)} $\g_4=\gl(1,\Co)\subset\u(1,1)$.

{\bf 5)} $\g_5=\gl(1,\Co)\zr\wedge^2\Real^2\subset\u(1,1)$.

{\bf 6)} $\g_6=\f\zr \wedge^2\Real^2$, where
$$\f=\left\{\left.\left(\begin{array}{cc}
a&c\\0&b\end{array}\right)\right|a,b,c\in\Real\right\}\subset\gl(2,\Real).$$

\end{theorem}

{\bf Proof.} The holonomy algebras of pseudo-Riemannian manifolds of
signature $(2,2)$ were classified in \cite{BBnn,Gh-Th}, see also
\cite{IRMA}. Comparing this result with Theorem \ref{th gener}, we
see that weakly irreducible not irreducible Berger algebras
$\g\subset\so(2,2)$ of Einstein type are exhausted by the algebras
form Theorem \ref{th 22} (note that the algebra $\f$ is conjugated
to the algebra $\g_1$, by this reason it does not appear in the
theorem). In order to complete the proof of the theorem we construct
Einstein not-Ricci flat metrics realizing the algebras from the
statement of the theorem as the holonomy algebras.

Consider on $\Real^4$ the coordinates $v_1,v_2,u_1,u_2$ and the
Walker metric
$$g=2dv_1du_1+2dv_2du_2+H_1(du_1)^2+2H_3du_1du_2+H_2(du_2)^2,$$
where $H_1$, $H_2$ and $H_3$ are functions. Let $\Lambda$ be a
non-zero real number. We claim that the following choice of the
functions provides Einstein metrics with the cosmological constant
$\Lambda$ and  with the corresponding holonomy algebras from the
statement of the theorem:

{\bf 1)} $H_1=\Lambda v_1^2+u_2^2$, $H_2=\Lambda v_2^2$, $H_3=0$;

\vskip0.2cm

{\bf 2)} $H_1=\frac{2\Lambda}{3} v_1^2,$ $H_2=\frac{2\Lambda}{3}
v_2^2$, $H_3=\frac{2\Lambda}{3} v_1v_2$;

\vskip0.2cm

{\bf 3)}  $H_1=\frac{2\Lambda}{3} v_1^2+v_1u_1u_2^2,$
$H_2=\frac{2\Lambda}{3} v_2^2-v_2u_1^2u_2$, $H_3=\frac{2\Lambda}{3}
v_1v_2$;

\vskip0.2cm

{\bf 4)} $H_1=\frac{\Lambda}{2}( v_1^2-v_2^2),$ $H_2=-H_1$,
$H_3=\Lambda v_1v_2$ ;

\vskip0.2cm

{\bf 5)} $H_1=\frac{\Lambda}{2}( v_1^2-v_2^2),$ $H_2=-H_1$,
$H_3=\Lambda v_1v_2+u_1u_2$;

\vskip0.2cm

{\bf 6)} $H_1=\Lambda v_1^2+{\rm ArcTan}(v_2)$, $H_2=\Lambda
v_2^2+\Lambda$, $H_3=0$.

\vskip0.2cm

We use a computer program to compute the Christoffel symbols, the
curvature tensor and the Ricci tensor of these metrics. We check
that the obtained metrics are Einstein with the cosmological
constant $\Lambda$. Let us explain, how to check that the holonomy
algebras of the metrics are the corresponding algebras from the
statement of the theorem. Since each metric under the consideration
is analytical, its  holonomy algebra $\g\subset\so(2,2)$  at the
point $0$ is generated by the elements of the form
$$\nabla^r_{{X_r}; \cdots ;X_1} R_0(X,Y),$$ where $r\geq 0$, and
$X$, $Y$, $X_1$,...,$X_r$ are tangent vector at the point $0$.

{\bf 1)} The vector fields $\partial_{v_1}$ and $\partial_{v_2}$ are
recurrent, hence $\g$ is contained in $\g_1$. The form of the
curvature tensor shows that the dimension of $\g$ is at least three.
This shows that $\g=\g_1$.

{\bf 2)}  The elements $R_0(X,Y)$, where $X$ and $Y$ are tangent
vectors at $0$, span $\g_2$, and it holds $\nabla R=0$.

{\bf 3)} The subalgebra $\g_3\subset\so(2,2)$ it the maximal
subalgebra preserving a fixed two-dimensional isotropic subspace and
its dimension is $5$. The elements of the form $R_0(X,Y)$ and
$\nabla_ZR_0(X,Y)$ generate a space of dimension $5$, consequently
$\g=\g_3$.

{\bf 4)} The considerations are the same as in the case 2).

{\bf 5)} The elements $R_0(X,Y)$ span $\g_5$.  Denote the
coordinates on $\Real^4$ by $x_1,\dots, x_4$. For each $i$, $1\leq
i\leq 4$, consider the matrix $\Gamma_i=(\Gamma_{ij}^k)_{j,k=1}^4$.
The values and all partial derivatives at $0$ of
$R(\partial_i,\partial_j)$ and $\Gamma_k$ belong to $\g_5$. Since
the covariant derivatives of $R$ are expressed in terms of the
partial derivatives and the Lie brackets of
$R(\partial_i,\partial_j)$ and $\Gamma_k$, it holds $\g=\g_5$.

{\bf 6)} The vector field $\partial_{v_1}$ is recurrent.  The
subalgebra $\g_6\subset\so(2,2)$ it the maximal subalgebra
preserving a fixed two-dimensional isotropic subspace and a line in
this subspace; its dimension is $4$. The elements of the form
$R_0(X,Y)$ and $\nabla_ZR_0(X,Y)$ generate a space of dimension $4$,
consequently $\g=\g_6$.

The theorem is proved. $\Box$

\vskip0.3cm

Consider now the general case. We use the notation of
Theorem~\ref{th gener}.

\begin{theorem}\label{th m=2} Weakly irreducible not irreducible
holonomy algebras $\g\subset\so(2,n+2)$  of Einstein not Ricci-flat
pseudo-Riemannian manifolds of signature $(2,n+2)$, $n\geq 1$, are
exhausted by the following algebras:

{\bf a)} $\g_a=\gl(1,\Real)\oplus\h\zr
\Real^{1,n+1}\subset\so(2,n+2)_\Real$, where $\h\subset\so(1,n+1)$
is the holonomy algebra of an Einstein not Ricci-flat Lorentzian
manifold that is the direct sum of an $\so(1,l+1)$ and a holonomy
algebra of an Einstein not Ricci-flat Riemannian manifold;

{\bf b)}
$\g_b=\gl(2,\Real)\oplus\h\zr(\Real^2\otimes\Real^n\zr\wedge^2\Real^2)\subset\so(2,n+2)_{\Real^2}$,
where $\h\subset\so(n)$ is the holonomy algebra of an Einstein not
Ricci-flat Riemannian manifold;

{\bf c)} $\g_c=\gl(1,\Co)\oplus\h\zr(\oplus_{\alpha=1}^t\N_\alpha\zr
\wedge^2\Real^2)\subset\so(2,n+2)_{\Real^2}$, where
$\h\subset\so(n)$ is the holonomy algebra of an Einstein not
Ricci-flat Riemannian manifold; $\N_\alpha\subset \Real^2\otimes
L_\alpha$ is a non-trivial subspace; if $\N_\alpha\neq
\Real^2\otimes L_\alpha$, then $\h_\alpha\subset\u(L_\alpha)$ and
either $\N_\alpha\simeq L_\alpha$, or $\N_\alpha\simeq\bar
L_\alpha$;

{\bf d)}
$\g_d=\gl(1,\Real)\oplus\gl(1,\Real)\oplus\h\zr(\oplus_{\alpha=1}^t\N_{1\alpha}\,\oplus\,
\oplus_{\alpha=1}^t\N_{2\alpha}\zr
\wedge^2\Real^2)\subset\so(2,n+2)_{\Real^2}$, where
$\h\subset\so(n)$ is the holonomy algebra of an Einstein not
Ricci-flat Riemannian manifold; each
$\N_{1\alpha},\N_{2\alpha}\subset L_\alpha$ is either $L_\alpha$ or
$0$; for each $\alpha$, at least one of $\N_{1\alpha}$ and
$\N_{2\alpha}$ equals to $L_\alpha$;

{\bf e)}
$\g_e=\f\oplus\h\zr(\oplus_{\alpha=1}^t\N_{1\alpha}\,\oplus\,
\oplus_{\alpha=1}^t\N_{2\alpha}\zr
\wedge^2\Real^2)\subset\so(2,n+2)_{\Real^2}$, where
$$\f=\left\{\left.\left(\begin{array}{cc}
a&c\\0&b\end{array}\right)\right|a,b,c\in\Real\right\},$$
 $\h\subset\so(n)$ is
the holonomy algebra of an Einstein  not Ricci-flat Riemannian
manifold; for each $\alpha$, $\N_{1\alpha}=L_\alpha$ and
$\N_{2\alpha}$ is either $L_\alpha$ or $0$.

\end{theorem}

{\bf Proof.}

{\bf 1. Algebraic classification.} Suppose that the maximal
dimension of a totally isotropic subspace preserved by
$\g\subset\so(2,n+2)$ is one, then from Theorem \ref{th gener} it
follows that $\g$ is the first algebra from the statement of the
theorem.

In other cases, $\g$ preserves a totally isotropic subspace of
dimension 2. Suppose that $\f\subset\gl(2,\Real)$ is irreducible.
Then either $\f=\gl(2,\Real)$, or $\g=\gl(1,\Co)$.

Suppose that $\g=\gl(1,\Co)$. By Theorem \ref{th gener},
$\N=\oplus_{\alpha=1}^t\N_\alpha$, and each
$\N_\alpha\subset\Real^2\otimes L_\alpha$ is a non-trivial
$\gl(1,\Co)\oplus\h_\alpha$-submodule. The space $\Real^2\otimes
L_\alpha$ can be considered as the complexification $\Co\otimes
L_\alpha$ of $L_\alpha$. This
$\gl(1,\Co)\oplus\h_\alpha=\Co\oplus\h_\alpha$-module is irreducible
if and only if  the complexified representation
$\Co\otimes\h_\alpha\subset\gl(\Co\otimes L_\alpha)$ is irreducible.
This is equivalent to the  non-existence of a complex structure on
$L_\alpha$. In the opposite case, $\Co\otimes L_\alpha\simeq
L_\alpha\oplus\bar L_\alpha$. We obtain the third algebra.

If $\g=\gl(2,\Real)$, then each $\gl(2,\Real)\oplus\h_\alpha$-module
$\Real^2\otimes L_\alpha$ is irreducible. . We obtain the second
algebra.

Suppose that $\f\subset\gl(2,\Real)$ is not irreducible. Then
either $\f=\gl(1,\Real)\oplus\gl(1,\Real)$, or
$\f=\gl(1,\Real)\oplus\gl(1,\Real)\oplus\f_{12}$, where
$\f_{12}=\Real$.

From Theorem \ref{th gener} it follows that
$\g\cap\N=\oplus_{\alpha=1}^t\N_{1\alpha}\,\oplus\,
\oplus_{\alpha=1}^t\N_{2\alpha}$, where each $\N_{i\alpha}$ is
either 0, or $L_\alpha$. Next, at least one of $\N_{1\alpha}$ or
$\N_{2\alpha}$ is non-zero, i.e. it equals to $L_\alpha$. Suppose
that $\N_{1\alpha}=\N_{2\alpha}=L_\alpha$ for some $\alpha$, then
$[\N_{1\alpha},\N_{2\alpha}]=\C=\wedge^2\Real^2$.

Suppose that $\g\cap\C=0$, then, for each $\alpha$, either
$\N_{1\alpha}=0$, or $\N_{2\alpha}=0$. Consequently $\g$ preserves
the non-degenerate subspaces $\Real p_1\oplus\,
\oplus_{\N_{2\alpha}=0}L_\alpha\,\oplus\Real q_1$, and $\Real
p_2\oplus\, \oplus_{\N_{1\alpha}=0}L_\alpha\,\oplus\Real q_2$. This
shows that $$\g=\gl(1,\Real)\oplus
(\oplus_{\N_{2\alpha}=0}\h_\alpha\zr L_\alpha)\oplus
\gl(1,\Real)\oplus (\oplus_{\N_{1\alpha}=0}\h_\alpha\zr L_\alpha)$$
is $\mathbb{Z}$-graded of depth 1.

If $\g\cap\C\neq 0$, then $\g$ is $\mathbb{Z}$-graded of depth 2,
hence it can not preserve a non-degenerate subspace of Lorentzian
signature. Clearly $\g$ cannot preserve any Euclidean subspace.
Thus, $\g$ is weakly irreducible. We obtain the algebra d).

Let $\f=\gl(1,\Real)\oplus\gl(1,\Real)\oplus\f_{12}$. Since
$[\f_{12},\N_{2\alpha}]=\N_{1\alpha}$, we see that
$\N_{1,\alpha}=L_\alpha$ for all $\alpha$. Suppose that
$\g\cap\C=0$. Then $\N_{2,\alpha}=0$ for all $\alpha$. In this case
$\g$ preserves the isotropic subspace spanned by the vectors $p_1$
and $q_2$. If we consider the vectors $p_1,q_2,q_1,p_2$ instead of
$p_1,p_2,q_1,q_2$, the we get the algebra d) ($\f_{12}$ becomes
$\C$).

If $\g\cap\C\neq 0$, then we get the algebra e). The fact that it is
weakly irreducible follows from the previous case.

{\bf 2. Construction of metrics.} To complete the proof of Theorem
\ref{th m=2}, we should show that each of the algebras
$\g\subset\so(2,n+2)$ from the statement of the theorem may appear
as the holonomy algebra of an Einstein not Ricci-flat
pseudo-Riemannian manifold. The metrics are constructed using the
idea from Section \ref{sec Lor}.

{\bf a)} For the case of the first algebra, the metric can be taken
as in Section \ref{sec Lor} with $h$ to be an Einstein not
Ricci-flat Lorentzian metric with the holonomy algebra
$\h\subset\so(1,n+1)$.

Next, let $h$ be an Einstein not Ricci-flat  Riemannian metric
defined on $\Real^n$ with the holonomy algebra $\h\subset\so(n)$ and
cosmological constant $\Lambda$. The de~Rham decomposition implies
that the coordinates can be divided into the groups
$$(x^1,...,x^n)=(x_1^\alpha,...,x^\alpha_ {n_\alpha})_{\alpha=1}^t$$
corresponding to the decomposition of $\h$ into irreducible parts.
Consider the coordinates $v_1,v_2,x_1,...,x_n,u_1,u_2$ on
$\Real^{n+4}$.

The algebras $\g_b,\g_{c},\g_d,\g_e$ may be obtained respectively
from the algebras $\g_3,\g_4,\g_1,\g_6$ by adding certain subspace
$\N\subset\Real^2\otimes\Real^n$. Fix one of the algebras from
Theorem \ref{th 22} and denote it be $\f$.  Denote by $f$ be the
corresponding metric constructed above. The holonomy algebra of the
product metric $f+h$ is $\f\oplus\h$. We will show that it is
possible to  twist slightly this metric in order to get an Einstein
not Ricci-flat metric with the corresponding holonomy algebra from
the statement of the theorem.

 Consider the metric
$$g=f+h+F_1(du^1)^2+2F_3du^1du^2+F_2(du^2)^2,$$
where
$$F_i=\sum_{\alpha=1}^{t} F_{i\alpha}(x^\alpha_1,...,x^\alpha_{n_\alpha}),\quad i=1,2,3,$$
are some function that should be chosen in a proper way in order to
add to the initial holonomy algebra $\f\oplus\h$ a  subspace
$\N_\alpha\subset \Real^2\otimes L_\alpha$. Since the coordinates
$x_1,...,x_n$ are separated into groups, it is enough to find a
construction for $t=1$.

The non-zero Chrystoffel symbols of the metric $f+h$ are some of
$$\Gamma^{v_k}_{v_iu_j},\quad \Gamma^{v_k}_{u_iu_j},\quad
\Gamma^{u_k}_{u_iu_j},\quad \Gamma^{x_a}_{x_bx_c}.$$

 The non-zero Chrystoffel symbols of the metric $g$ are the following non-zero
Chrystoffel symbols of the metric $f+h$:
$$\Gamma^{v_k}_{v_iu_j},\quad
\Gamma^{u_k}_{u_iu_j},\quad \Gamma^{x_a}_{x_bx_c},$$ some of the
modified $\Gamma^{v_k}_{u_iu_j}$,   and some of the new symbols
$$\Gamma^{v_k}_{x_au_j},\quad \Gamma^{x_a}_{u_iu_j}.$$
In particular, it holds
$$\Gamma^{v_i}_{x_au_j}=\frac{1}{2}\partial_{x_a}F_{ij},$$
where we use the denotation $F_{11}=F_1,$ $F_{12}=F_{21}=F_3$,
$F_{22}=F_2$. Since we are interested in the projection of the
holonomy algebra onto $\Real^2\otimes\Real^n$, the symbols
$\Gamma^{v_i}_{x_au_j}$ are the only  we should carry  about.
Similarly, we are interested in the components
$\nabla_{\bullet;\cdots;\bullet}R^{v_i}_{x_ax_bu_j}$ of the
curvature tensor and its covariant derivatives. It holds
$$R^{v_i}_{x_ax_bu_j}=\frac{1}{2}\partial_{x_a}\partial_{x_b}F_{ij}-\frac{1}{2}\sum_{c=1}^n\Gamma^{x_c}_{x_ax_b}\partial_{x_c}F_{ij}.$$
Note that $R^{v_i}_{x_au_ku_j}=0$.

Consider now the algebras from the statement of the theorem.

{\bf b)} We take $F_2=F_3=0$. The Einstein equation (with the
cosmological constant $\Lambda$) is equivalent to the condition
$$\Delta_hF_1=0,$$ where $\Delta_h$ is the Laplace operator with
respect to the metric $h$. If we take the function $F_1$
sufficiently general (such that $R^{v_i}_{x_ax_bu_j}\neq 0$ for some
indices), then the holonomy algebra $\g$ of that metric at the point
0 contains a non-trivial intersection with $\Real^2\otimes\Real^n$.
Since the holonomy algebra $\g$ contains also $\gl(2,\Real)\zr
\wedge^2\Real$ and $\h$, by the proof of the first part of the
theorem, $\g=\g_b$.

{\bf c)} If $\N=\Real^2\otimes\Real^n$, then the construction is the
same as in the previous case. Suppose that $n=2l$, $\h\subset\u(l)$,
and $\N=\Co^l\subset \Real^2\otimes \Real^{2l}$. More precisely,
$\N$ consists of the matrices
$$\left(\begin{matrix}X&-Y\\Y&X\end{matrix}\right),$$ where
$X,Y\in\Real^l$. Denote by $J$ the K\"ahlerian structure of the
metric $h$. Let the coordinates $x_1,\dots x_{2l}$ satisfy
$$J\partial_{x_a}=\partial_{x_{a+l}},\quad J\partial_{x_{a+l}}=-\partial_{x_a},\quad a=1,\dots,l.$$
Let $F_2=-F_1$ and
$$F_1(x_1,\dots,x_{2l})=F_3(x_{l+1},\dots,x_{2l},-x_1,\dots,
-x_l).$$ The Einstein equation is equivalent to the condition
$\Delta_hF_3=0$. We get
$$\Gamma^{v_1}_{x_au_j}=\Gamma^{v_2}_{x_{a+l}u_j},\quad \Gamma^{v_2}_{x_au_j}=-\Gamma^{v_1}_{x_{a+l}u_j},\quad a=1,\dots,l,\quad j=1,2.$$
This implies
$$R^{v_1}_{x_ax_bu_j}=R^{v_2}_{x_{a+l}x_bu_j},\quad R^{v_2}_{x_ax_bu_j}=-R^{v_1}_{x_{a+l}x_bu_j},\quad a,b=1,\dots,l,\quad j=1,2.$$
Consequently the same  condition satisfy all covariant derivatives
of the curvature tensor at the point $0$, this implies that
$\g=\g_c$.

{\bf d,e)} Let $F_3=0$. The Einstein equation is equivalent to the
condition $\Delta_hF_1=\Delta_hF_2=0$. If we take both $F_1$ and
$F_2$ sufficiently general, then $\N=\Real^2\otimes \Real^{n}$. If
we take $F_2=0$ and $F_1$ sufficiently general, then
$\N=\N_{11}=\Real^n\subset\Real^2\otimes \Real^{n}$.

The theorem is proved.
 $\Box$

\section{Para-quaternionic-K\"ahlerian manifolds}\label{paraqK}
Quaternionic-K\"ahlerian manifolds of non-zero scalar curvature
are always Einstein and not Ricci-flat. Holonomy algebras of these
manifolds in arbitrary signature are classified in \cite{Quat}.
This result may be also deduced from Theorem \ref{th gener} in the
same way as we do now for the case of
para-quaternionic-K\"ahlerian manifolds.

Recall that a para-quaternionic-K\"ahlerian manifold \cite{A-C,DJS}
is a pseudo-Riemannian manifold $(M,g)$ that admits a parallel
three-dimensional subbundle of the bundle of endomorphisms of the
tanget bundle of $M$ locally spanned by endomorphisms $I$, $S$, $T$
preserving $g$ and satisfying the relations of the split quaternions
$$-I^2=S^2=T^2=\id,\quad IS=T=-SI.$$ Note that the endomorphisms
$I$, $S$, $T$ generate the Lie algebra isomorphic to $\sl(2,\Real)$.

Equivalently, a para-quaternionic-K\"ahlerian manifold is a
pseudo-Riemannian manifold with the holonomy group contained in
$${\rm Sp}(2n,\Real)\cdot {\rm SL}(2,\Real)\subset{\rm
SO}(2n,2n).$$ For the holonomy algebra we get
$$\g\subset\sp(2n,\Real)\oplus\sl(2,\Real)\subset\so(2n,2n).$$ In \cite{A-C} it is shown that each
para-quaternionic-K\"ahlerian manifold of non-zero scalar
curvature is Einstein, locally indecomposable (i.e. with weakly
irreducible holonomy algebra), and its holonomy algebra contains
the subalgebra $\sl(2,\Real)$.

Let $(M,g)$ be a para-quaternionic-K\"ahlerian manifold of
non-zero scalar curvature and let
$\g\subset\sp(2n,\Real)\oplus\sl(2,\Real)$ be its holonomy
algebra. Suppose that $\g$ is not irreducible, i.e. it preserves
an isotropic subspace of the tangent space. Since $\g$ contains
the subalgebra $\sl(2,\Real)$, which does not annihilate any
subspace of the tangent space, Theorem \ref{th gener} implies that
$r+s=0$, $m=2n$, $\k=\k_0=\k_1\subset\gl(2n,\Real)$ is
irreducible, and
$$\g\subset\gl(2n,\Real)\zr\wedge^2\Real^{2n}.$$
Since $\g$ stabilizes $\sl(2,\Real)$, i.e.
$[\g,\sl(2,\Real)]\subset \sl(2,\Real)$, we get that
$$\pr_{\gl(2n,\Real)}\g\subset\gl(n,\Real)\oplus\sl(2,\Real),$$
where the representation of $\gl(n,\Real)\oplus\sl(2,\Real)$ in
$\Real^{2n}$ is given by the tensor products of representations:
$\Real^{2n}=\Real^n\otimes\Real^2$. Next,
$$\wedge^2\Real^{2n}=(\wedge^2 \Real^n\otimes\odot \Real^2)\oplus (\odot^2 \Real^n\otimes\wedge^2 \Real^2)$$
as the $\gl(n,\Real)\oplus\sl(2,\Real)$-module. Note that the
$\sl(2,\Real)$-module $\wedge^2 \Real^2$ is trivial.  Since $\g$
stabilizes $\sl(2,\Real)$, we get
$$\pr_{\wedge^2\Real^{2n}}\g\subset (\odot^2 \Real^n\otimes\wedge^2
\Real^2).$$ Thus,
$$\g=\k\zr\C,\quad \k\subset \gl(n,\Real)\oplus\sl(2,\Real),\quad
\C\subset \odot^2 \Real^n\otimes\wedge^2 \Real^2.$$ Let
$\g_0=\gl(n,\Real)\oplus\sl(2,\Real)\subset\gl(2n,\Real)$. The
corresponding subalgebra $\g_0\subset\so(2n,2n)$ is the holonomy
algebra of the para-quaternionic-K\"ahlerian symmetric space
$${\rm SL}(n+2,\Real)/({\rm GL}(n,\Real)\cdot{\rm SL}(2,\Real)).$$
Let $R_0\in\R_1(\g_0\subset\so(2n,2n))$ be the corresponding
algebraic curvature tensor. From Lemma \ref{lem Wu g nn} it
follows that $$\R(\gl(n,\Real)\subset\so(2n,2n))=0.$$ This implies
that $$\R(\g_0\subset\so(2n,2n))=\Real R_0.$$ This and the proof
of Lemma \ref{lem vid f} show that $\k=\g_0$. Finally note that
the $\g_0$-module $\odot^2 \Real^n\otimes\wedge^2 \Real^2$ is
irreducible, consequently, either $\C=0$, or $\C=\odot^2
\Real^n\otimes\wedge^2 \Real^2$. We have proved the following
theorem.

\begin{theorem} Let $(M,g)$ be a para-quaternionic-K\"ahlerian
manifold  of non-zero scalar curvature and dimension $4n$. If its
holonomy algebra $\g$ is not irreducible, then it preserves an
isotropic subspace of dimension $2n$, and it coincides with one of
the following subalgebras of $\so(2n,2n)_{\Real^{2n}}$:

$\bullet$ $\gl(n,\Real)\oplus\sl(2,\Real)$,

$\bullet$ $(\gl(n,\Real)\oplus\sl(2,\Real))\zr (\odot^2
\Real^n\otimes\wedge^2 \Real^2)$.
\end{theorem}

Let us note that the irreducible holonomy algebras of
 para-quaternionic-K\"ahlerian
manifolds  of non-zero scalar curvature are exhausted by
$\sp(2n,\Real)\oplus\sl(2,\Real)$ and by the holonomy algebras of
para-quaternionic-K\"ahlerian symmetric spaces; the list of
para-quaternionic-K\"ahlerian symmetric spaces may be found in
\cite{A-C,DJS}, or it may be deduced from \cite{Ber57}; the only
symmetric space with non-irreducible holonomy algebra is ${\rm
SL}(n+2,\Real)/({\rm GL}(n,\Real)\cdot{\rm SL}(2,\Real)).$

\section{Final remarks}\label{sec Rem}

It is an open problem to  construct  examples of Einstein not
Ricci-flat  metrics with  the holonomy algebras
$\g\subset\so(m+r,m+s)$ from Theorem \ref{th gener} for signatures
different from $(1,N)$ and $(2,N)$. The following idea generalizing
the construction from Section \ref{sec m=2} may be used. Let
$\g\subset\so(m+r,m+s)$ be an algebra from Theorem \ref{th gener}.
Consider the coordinates
$$(v_1^i,\dots,v_{m_i}^i)_{i=1}^k, (x_1^\alpha,\dots
x_{n_\alpha}^\alpha)_{\alpha=1}^t,(u_1^i,\dots,u_{m_i}^i)_{i=1}^k$$
on $\Real^{m+r,m+s}$. Let
$$f_i=\sum_{a=1}^{m_i}2dv^i_adu^i_a+\sum_{a,b=1}^{m_i}H_{ab}(v^i_1,...,v^i_{m_i})du^i_adu^i_b$$
be an Einstein metric with the cosmological constant $\Lambda$ and
the holonomy algebra $\f_i\subset\so(m_i,m_i)$. Let
$$h_\alpha=\sum_{a,b=1}^{n_\alpha}h^\alpha_{ab}(x^\alpha_1,...,x^\alpha_{n_\alpha})dx^\alpha_adx^\alpha_b$$
be an Einstein metric with the cosmological constant $\Lambda$ and
the
 holonomy algebra $\h_\alpha\subset\so(n_\alpha)$. Such metrics
 exist, one may consider e.g. metrics of symmetric spaces.
The metric $$g=f+h=\sum_{i=1}^kf_i+\sum_{\alpha=1}^t h_\alpha$$ is
Einstein with the holonomy algebra $\f_0\oplus\h$. This metric is
decomposable. Next, if some of $\f_{ij}$, $N_{i\alpha}$ or $\C_{ij}$
is non-trivial, then we farther curve the metric adding to $g$ the
following terms:
$$f_{ij}=\sum_{a,b=1}^{m_j}F_{ij}^{ab}(v^i_1,...,v^i_{m_i})du^i_adu^j_b,\quad
n_{i\alpha}=\sum_{a,b=1}^{m_i}
N_{i\alpha}^{ab}(x^\alpha_1,...,x^\alpha_{n_\alpha})du^i_adu^i_b,$$
$$c_{ij}=\sum_{a=1}^{m_i}\sum_{b=1}^{m_j}C_{ij}^{ab}(u^i_1,...,u^i_{m_i},u^j_1,...,u^j_{m_j})du^i_adu^j_b.$$
The additional terms will give some additional components to the
curvature tensor of the metric $f+h$, these components will take
values in $V_j^*\otimes V_i$, $V_i\otimes L_\alpha$, and $V_i\otimes
V_j$, respectively.  The Einstein condition and the condition on the
holonomy algebra to coincide with $\g$ will give equations on the
functions $F_{ij}^{ab}$,  $N_{i\alpha}^{ab}$, $C_{ij}^{ab}$. Above
we have seen that for each of $F_{ij}^{ab}$,  $N_{i\alpha}^{ab}$,
$C_{ij}^{ab}$, $i\neq j$, there are only three possibilities: it can
be trivial, coincide with the one of the corresponding spaces
$V_j^*\otimes V_i$, $V_i\otimes L_\alpha$,  $V_i\otimes V_j$ or with
one of the corresponding spaces $V_j^*\otimes_\Co V_i$,
$V_i\otimes_\Co L_\alpha$, $V_i\otimes_\Co V_j$. In the first case,
the corresponding functions may be chosen to be zero; in the second
case, the corresponding functions may be chosen to be harmonic and
sufficiently general; in the third case, the corresponding functions
may be found in the similar way as in case c) from  the proof of
Theorem \ref{th m=2}. A complication may appear, when one considers
the $\f_i$-modules $\C_{ii}$. Each $\C_{ii}$ is a submodule of
$\wedge^2V_i$. One should go through the list of the representations
$\f_i\subset\gl(m_i,\Real)$ from Section \ref{sec nn}, and in each
case describe all submodules of the wedge product $\wedge^2V_i$. It
is a challenge to find proper functions $C_{ii}^{ab}$ for each of
these submodules. One example of the module $\C_{ii}$ provides
Section \ref{paraqK}. Since the rigorous construction and its
justification require complicated technical work, they should be
done in a separate paper.

We complete the paper by giving two remarks. First, the holonomy
algebra of an Einstein not Ricci-flat pseudo-Riemannian  symmetric
space is reductive \cite{AlIRMA}. Thus the only holonomy algebras of
symmetric spaces appearing in Theorem \ref{th gener} are contained
in $\so(n,n)$ and have been discussed in Section \ref{sec nn}.
Finally, note that in \cite{GalmNotes} we show that in the general
non-Einstein case the subalgebra $\h\subset\so(r,s)$ associated to a
holonomy algebra $\g\subset\so(r+m,s+m)$ may be absolutely
arbitrary, i.e. the Einstein case is in sharp contrast to the
non-Einstein case.

\vskip0.3cm

{\bf Acknowledgements.} The author is thankful to Lionel
B\'erard~Bergery for sending him the paper in preparation
\cite{BBnn}. The work is partially supported by the grant  no.
18-00496S of the Czech Science Foundation.


\affiliationone{University of Hradec Kr\'alov\'e\\ Faculty of
Science\\ Rokitansk\'eho 62\\ 500~03 Hradec Kr\'alov\'e\\  Czech
Republic
   \email{anton.galaev(at)uhk.cz}}

\end{document}